\documentclass{article}

\usepackage[tbtags]{amsmath}
\usepackage{amsfonts,amssymb}
\usepackage{theorem}

\begin{document}

\title{Number of binomial coefficients divided by a fixed power of a prime}

\author{William B.~Everett\\
Chernogolovka, Moscow Oblast, Russia}

\maketitle

\begin{abstract}
We state a general formula for the number of binomial coefficients $n$ choose
$k$ that are divided by a fixed power of a prime $p$, i.e., the number of
binomial coefficients divided by $p^j$ and not divided by $p^{j+1}$.
\end{abstract}

Let $n$ be a natural number and $p$ be a prime. Let $\theta_j(n)$ denote the
number of binomial coefficients
$_nC_k=\left(\begin{smallmatrix}n\\k\end{smallmatrix}\right)$ such that $p^j$
divides $_nC_k$ and $p^{j+1}$ does not divide $_nC_k$.

We represent $n$ in the base $p$: $n=c_0+c_1p+c_2p^2+\dots+c_rp^r$,
$0\le c_i<p$, $i=0,1,\dots,r$, $c_r\ne0$.

Let $W$ be the set of $r$-bit binary words, i.e.,
\[
W=\bigl\{\mathbf{w}=w_1w_2\dots w_r\colon w_i\in\{0,1\},\;1\le i\le r\bigr\}.
\]
We partition $W$ into $r{+}1$ subsets $W_j$, $0\le j\le r$, where
\[
W_j=\biggl\{\mathbf{w}\in W\colon \sum_{i=1}^rw_i=j\biggr\}.
\]

We define the functions $F(\mathbf{w})$, $L(\mathbf{w})$, and
$M(\mathbf{w},i)$ as follows:
\begin{align*}
&F(\mathbf{w})=\begin{cases}
c_0+1&\text{if }w_1=0,\\
p-c_0-1&\text{if }w_1=1,\end{cases}
\\
&L(\mathbf{w})=\begin{cases}
c_r+1&\text{if }w_r=0,\\
c_r&\text{if }w_r=1,\end{cases}
\\
&M(\mathbf{w},i)=\begin{cases}
c_i+1&\text{if }w_i=0\text{ and }w_{i+1}=0,\\
p-c_i-1&\text{if }w_i=0\text{ and }w_{i+1}=1,\\
c_i&\text{if }w_i=1\text{ and }w_{i+1}=0,\\
p-c_i&\text{if }w_i=1\text{ and }w_{i+1}=1.\end{cases}
\end{align*}

The general formula for $\theta_j(n)$ is
\begin{equation}
\theta_j(n)=\sum_{\mathbf{w}\in W_j}F(\mathbf{w})L(\mathbf{w})
\prod_{i=1}^{r-1}M(\mathbf{w},i).
\label{eq1}
\end{equation}
Obviously, we have the sum of $_rC_j$ terms and each term is the product of
$r{+}1$ factors. It is easy to establish that if $p-c_\ell-1$=0 for some
$\ell$, then $p-c_i-1=0$ for all $i\le\ell$ (and some terms may vanish from
the sum). It is also easy to establish that if $p-c_\ell-1=1$ for some
$\ell$, then $p-c_i-1=1$ for all $i\le\ell$ (and the number of contributing
factors in some terms is reduced). This means that the formula can be
simplified for $n$ of certain special forms.

Formula~\eqref{eq1} reproduces known formulas for some particular values. For
example, for $j=0$, we obtain the known formula~\cite{1}
\begin{align*}
\theta_0(n)&=(c_0+1)(c_r+1)(c_1+1)\cdots(c_{r-1}+1)
\\
&=(c_0+1)(c_1+1)\cdots(c_r+1).
\end{align*}
For $j=1$, we obtain the known formula~\cite{2}
\begin{align*}
\theta_1(n)={}&(c_0+1)c_r(c_1+1)(c_2+1)\cdots(c_{r-2}+1)(p-c_{r-1}-1)
\\
&{}+(c_0+1)(c_r+1)(c_1+1)(c_2+1)\cdots(p-c_{r-2}-1)c_{r-1}
\\
&{}+\cdots
\\
&{}+(c_0+1)(c_r+1)(p-c_1-1)c_2(c_3+1)\cdots(c_{r-1}+1)
\\
&{}+(p-c_0-1)(c_r+1)c_1(c_2+1)\cdots(c_{r-1}+1)
\\
={}&\sum_{k=0}^{r-1}(c_0+1)\cdots(c_{k-1}+1)(p-c_k-1)c_{k+1}
(c+{k+2}+1)\cdots(c_r+1).
\end{align*}
Other particular formulas for $\theta_j(n)$ can be found in~\cite{3}
and~\cite{4}.

\end{document}